\documentclass{amsart}
\usepackage{amssymb}

 \newtheorem{theorem}{Theorem}
 \newtheorem{lemma}[theorem]{Lemma}
  \newtheorem{proposition}[theorem]{Proposition}
 \newtheorem{corollary}[theorem]{Corollary}
 \newtheorem{remark}[theorem]{Remark}
 \newtheorem{example}[theorem]{Example}
 \newtheorem{definition}[theorem]{Definition}
 \newtheorem{conjecture}[theorem]{Conjecture}
 
 \newtheorem{problem}[theorem]{Problem}
 
 \newtheorem{question}[theorem]{Question}

\newcommand{\be}{\begin}
\newcommand{\en}{\end}
\newcommand{\bpr}{\begin{proof}}
\newcommand{\epr}{\end{proof}}
\newcommand{\beq}{\begin{equation}}
\newcommand{\eeq}{\end{equation}}
\newcommand{\bthm}{\begin{theorem}}
\newcommand{\ethm}{\end{theorem}}
\newcommand{\blem}{\begin{lemma}}
\newcommand{\elem}{\end{lemma}}
\newcommand{\bpro}{\begin{proposition}}
\newcommand{\epro}{\end{proposition}}
\newcommand{\bcor}{\begin{corollary}}
\newcommand{\ecor}{\end{corollary}}
\newcommand{\brem}{\begin{remark}}
\newcommand{\erem}{\end{remark}}
\newcommand{\bexa}{\begin{example}}
\newcommand{\eexa}{\end{example}}
\newcommand{\bdf}{\begin{definition}}
\newcommand{\edf}{\end{definition}}
\newcommand{\bcon}{\begin{conjecture}}
\newcommand{\econ}{\end{conjecture}}
\newcommand{\bque}{\begin{question}}
\newcommand{\eque}{\end{question}}
\newcommand{\bprb}{\begin{problem}}
\newcommand{\eprb}{\end{problem}}

 \newcommand{\Z}{{\mathbb Z}}

\newcommand{\p}{\partial}

\newcommand{\la}{\left\langle}
\newcommand{\ra}{\right\rangle}

\title{Distributive Products and Their Homology}
\author{J\'ozef H. Przytycki, Adam S. Sikora}
\address{Dept. of Mathematics, GWU, Monroe Hall, Room 240, 2115 G Street NW,
Washington, DC 20052, USA and Department of Mathematics, Physics and Computer Science, Gda\'nsk University, Gda\'nsk, Poland}
\address{244 Math Bldg, SUNY Buffalo, Buffalo NY 14260, USA}


\thanks{The first author was partially supported by the NSA-AMS 091111 grant, by the
Polish Scientific Grant: Nr. N-N201387034, by the GWU REF grant, and by the CCAS/UFF award.}

\date{}

\begin{document}
\thispagestyle{empty}

\begin{abstract} We develop a theory of sets with distributive products (called shelves and multi-shelves) and of their homology. We relate the shelf homology to the rack and quandle homology.
\end{abstract}

\maketitle

%
\section{Binary distributive operations}
%

Consider binary operations $\star:X\times X\to X$ on a set $X$ with the composition
$$x\star_1\star_2y=(x\star_1y)\star_2y.$$
Since the composition operation is associative, these binary operations form a monoid, $M(X),$ with the identity $x\star_0 y=x.$

A collection of operations $\star_i$, $i\in I$ is called {\em mutually (right) distributive} if $$(a\star_ib)\star_jc=(a\star_jc)\star_i(b\star_jc)$$ for all $i,j\in I.$
A direct computation gives:

\blem If $\{\star_i\}_{i\in I}$ are mutually-distributive then the set of all compositions of $\star_i$'s is mutually distributive. Hence every maximal set of mutually-distributive products is a monoid.
\elem


We say that $\star$ is {\em (right) self-distributive} or, simply, {\em distributive}, if
\beq\label{sdist}
(x\star y)\star z=(x\star z)\star(y\star z),
\eeq
that is, if the one element set $\{\star\}$ is mutually distributive. In this case, $(X,\star)$ is called a {\em shelf}, \cite{Cr}. (Hence (\ref{sdist}) might be called ``shelf-distributivity".)
We abbreviate $x\star y$ to $xy$ whenever it does not lead to a confusion.
Following the ``left normed" convention, we will abbreviate the products
$$(...((x_0\star x_1)\star x_2)... \star x_{n-1})\star x_n$$
by $x_0\star x_1\star ...\star x_n.$

A rich source of examples is given by quandles (those arise for example from groups) and, more generally, racks. There is a rich literature on these structures. In this paper, however, we will be primarily interested in shelves which are not racks. Below we present several examples of shelves which are, in general, not racks:

\bexa\label{shelf-ex} Every set $X$ with the one of the following products is a shelf:\\
(1) $x\star y=f(x)$, for some $f: X\to X.$\\
(2) $x\star y=g(y)$ for some $g: X\to X$ such that $g^2=g.$\\
(3) Given $A\subset X$, $x\star y=\begin{cases} x & \text{if $y\in A$}\\  y & \text{otherwise.}\end{cases}$\\
(4) Let $b\in X$, $x\star y=\begin{cases} g(y) & \text{if $x=b$}\\  y & \text{otherwise,}\end{cases}$\\ for any $g: X\to X,$ such that $g^{-1}(b)={b}.$
\eexa

The last example generalizes as follows:\\
(4') If $X$ is a disjoint union of $A_i$'s for $i\in I$ and
$g_i:X\to X$ are functions for $i\in I$ such that $g_i(A_j)\subset A_j$ and $g_i=g_jg_i$ on $A_j$
for all $i,j\in I$ then $x\star y=g_i(y)$, for $i$ such that $x\in A_i$, is a distributive product on $X$.

\bpr If $x_1\in A_{i_1}, x_2\in A_{i_2}, x_3\in A_{i_3}$, then
$$(x_1\star x_2)\star x_3=g_{i_1}(x_2)\star x_3=g_{i_2}(x_3)$$ and
$$(x_1\star x_3)\star (x_2\star x_3)=g_{i_1}(x_3)\star g_{i_2}(x_3)=g_{i_3}g_{i_2}(x_3)$$ coincide.
\epr

\bexa\label{shelf-ex-omega} Let $X$ be a collection of subsets of $\Omega.$\\
(1) If $X$ is closed under intersections then $(X, \cap)$ is a shelf.
(More generally, every semi-lattice is a shelf.)\\
(2) If $X$ is closed under set subtractions then $X$ with operation of subtraction, $x*y=x-y,$ is a shelf.
\eexa

We call a collection of mutually distributive products on a set, $(X,\star_i, i\in I)$, {\em a multi-shelf}. Clearly, every self-distributive product is mutually distributive with the identity product. Here are two more interesting examples of multi-shelves:

\bexa\label{mshelf-ex}
(1) If a shelf $(X,\star)$ is a spindle (i.e. $x \star x=x$ for every $x\in X$) then $\star$
is mutually distributive with $x\star_\sim y=y$.\\
(2) Any subset $X$ of $2^\Omega$ closed under unions and intersections, together with operations $\cap, \cup,$ the ``identity" product, $x\star y=x,$ and $x\star_\sim y=y$ is a multi-shelf (which we call the Boolean multi-shelf). By Birkhoff's and Stone's theorems, \cite{Bi,St}, every distributive lattice is of the form $(X,\cap, \cup).$
\eexa

Finally, we have

\bdf The free shelf $X(S)$ on a set $S$ is the free magma, $F(S)$, on $S$, quotiented by the smallest equivalence relation relating two words in $F(S)$, if they are related by the distribution operation $(w_1w_2)w_3=(w_1w_3)(w_2w_3).$
\edf

\be{pro}
There are precisely $6$ different two element shelves up to an isomorphism.
Two of them are given by Example \ref{shelf-ex-omega}(1) and (2) for $X=\{\emptyset,\Omega\}.$
The other four are of the types described in Example \ref{shelf-ex}(1) and \ref{shelf-ex}(2) for $X=\{0,1\}$:
(a) $f=0$, (b) $f=Id_X$, (c) $f(x)=1-x$, (d) $g=Id_X$.
\en{pro}

\be{proof} by direct enumeration.
\en{proof}

We present below some simple ways of building ``bigger" shelves from smaller ones.
We leave the proofs of the following two propositions to the reader.

\be{pro}\label{s_retract}
Every self-distributive product on $A$ can be extended to a self-distributive product on any
set $X$ containing $A$ as follows:
Consider any map $r:X\to A$ which restricts to the identity on $A$
and set $x\star y=r(x)\star r(y)$ for all $x,y\in X.$
In this case we say that $r$ is a strong retract\footnote{More generally, $r:X\to A\subset X$ is a retract
if $r$ is the identity on $A$ and $r(x\star y)=r(x)\star r(y).$} of $X$ onto $A$.
\en{pro}

In particular, for the shelves in Example \ref{shelf-ex}(2) $g$ is a strong retraction onto $A=g(X)$
with the product on $A$ given by $x\star y=y.$

\be{pro} If $(X_i,\star_i)$ are shelves for $i\in I$ then\\
(1) their disjoint union $X=\coprod_{i\in I} X_i$ with
$$x\star x'=\begin{cases} x \star_i x' & \text{iff}\ x,x'\in X_i\ \text{for some $i$}\\
x & \text{iff}\ x\in X_i, x'\in X_j, i\ne j\\
\end{cases}$$
is a shelf as well.\\
(2) $X= \prod_{i\in I} X_i$ with
$$(x_i,\, i\in I)\star(x_i',\, i\in I)=(x_i \star_i x_i',\, i\in I)$$
is a shelf as well.
\en{pro}

\bprb
Find other "natural" families of shelves which are not racks.
\eprb


%
\section{The simplicial complex and the homology of a shelf}
%

\subsection{Shelf complex}
Given a shelf $X$ consider a simplicial complex $S(X)$ whose vertices are elements of $X$ and whose simplices are of the form
$$\sigma_{x_0,...,x_n}= \la x_0x_1...x_{n-1}x_n,x_1...x_{n-1}x_n,...,x_{n-1}x_n,x_n\ra,$$
for $x_0,...,x_n\in X$ such that $x_0x_1...x_{n-1}x_n,x_1...x_{n-1}x_n,...,x_{n-1}x_n,x_n$
are distinct. The distributivity of the product implies that faces of simplices are simplices as well. This is {\em the shelf complex of $X$}.
Since every simplex of a simplicial complex is determined by its vertices,
$\sigma_{x_0,...,x_n}$ may denote the same simplex for different $n$-tuples $(x_0,...,x_n).$ We will call the homology of that complex the {\em simplicial shelf homology of $X$}.

\bexa Let $(m_{ij})=\left[\begin{array}{cccc} 1 & 3 & 3 & 4\\
1&2&3&4\\ 1 & 3 & 3 & 4\\ 3 & 1&3&4\\
\end{array}\right]$.
For the shelf $X=\{1,2,3,4\}$ with multiplication given by $i*j=m_{ij},$
$S(X)$ is composed of the $1$-dimensional triangle with
vertices $1,2,3$ together with a disjoint vertex $4.$ (The edges are given by $\la 4\star 1,1 \ra,
\la 4\star 2,2 \ra,\la 1\star 2,2 \ra$. One can check that $\la 1,2,3\ra\not \in S(X)$. Since $x\star 4=4$ for all $x,$ the vertex $4$ is disconnected.)
Hence the first simplicial homology of $X$ is $\Z.$
\eexa

Alternatively, one can consider a CW-complex, $S'(X)$, whose $n$-cells, $\sigma_{x_0,...,x_n},$ are in $1$-$1$ correspondence with all $(n+1)$-tuples of elements of $X$.

\subsection{Shelf homology}
Here is another version of homology for shelves:

Let $C_n(X)$ be the free abelian group with basis
$\{(x_0,...,x_n): x_0,...,x_n\in X\}$ and let $d_n:C_n(X)\to C_{n-1}(X)$ be given by $$d_n(x_0,...,x_n)=\sum_{i=0}^n (-1)^i (x_0\star x_i,...,x_{i-1}\star x_i,x_{i+1},...,x_n).$$
(For example, $d_1(x_0,x_1)=(x_1)-(x_0x_1),$ $d_2(x_0,x_1,x_2)=(x_1,x_2)-(x_0x_1,x_2)+(x_0x_2,x_1x_2).$)
Additionally, we assume that $C_{-1}(X)=\Z$ and $d_0:C_0(X)=ZX\to C_{-1}(X)$ sends every $(x)$ to $1.$

\blem\label{dd=0}
$(C_*(X),d_*)$ is a chain complex.
\elem

\bpr
Let $$d_{n,i}(x_0,...,x_n)=(x_0\star x_i,...,x_{i-1}\star x_i,x_{i+1},...,x_n).$$
Then $d_n=\sum_{i=0}^n (-1)^i d_{n,i}$ and it is easy to see that
$d_{n-1,i+1}d_{n,j}=d_{n-1,j}d_{n,i},$ for $i\geq j$.
Hence, the $(n+1)n$ summands of $d_{n-1}d_n(x_0, ...,x_n)$ can be matched in pairs of coinciding terms of opposite signs.
\epr

We call the homology groups of this chain complex, the {\em (algebraic) shelf homology of $X$.} Its generalizations and relations to rack and quandle homologies will be explained in Section \ref{s_multi}.

Let $(C_*(S(X)),d_*)$ be the chain complex of the simplicial complex $S(X).$ We  use here the convention that $C_n(S(X))$ is the quotient of the free abelian group on all $n$-tuples of vertices $(x_1,...,x_n)$ of $S(X)$ by the relation $$(x_{\sigma(1)},...,x_{\sigma(n)})=sign(\sigma)(x_1,...,x_n),$$
for every permutation $\sigma\in S_n.$ (In particular, $(x_1,...,x_n)=0$ if $x_i=x_j$ for some $i\ne j.$)
Note that there is a chain map from $(C_*(X),d_*)$ to $(C_*(S(X)),d_*)$ sending $(x_0,...,x_n)$ to
$(x_0x_1...x_n, x_1...x_n,..., x_n).$ This map induces a homomorphism of corresponding homology groups.

\subsection{$H_0$ and Left Orbits}

Let $\sim$ be the smallest equivalence relation
on $X$ such that $x\sim y\star x$ for all $x,y\in X.$ We call the equivalence classes of this relation the {\em left orbits}. A shelf with a single left orbit is called {\em left connected}.

\be{rem} (1) The following shelves are left connected:\\
(a) Racks\\
(b) The shelves in Examples \ref{shelf-ex}(1),(3).\\
(c) The shelves in Example \ref{shelf-ex-omega}(1) and (2): (Proof for (1): $x\sim y\cap x=x\cap y\sim y$. Proof for (2): $x-x=\emptyset.$ Hence $x\sim \emptyset$ for every $x\in X.$)\\
(2) The left orbits in Example \ref{shelf-ex}(2) are in bijection with the elements of the image of $g$.
\en{rem}

\be{rem}\label{r_orb} Since $x'\sim x$ implies $y'\star x'\sim x'\sim x \sim y\star x$ for every $y,y'\in X$, the product on $X$ descends to a self-distributive product on the set $O$ of the left orbits in $X$ and the natural quotient map $\pi: X\to O$ is a homomorphism.
The induced product on $O$ is $[x]\star [y]=[y]$ (as in Example 2(2) for $g=Id$).
\en{rem}

\bpro There is a natural bijection between left orbits of $X$ and the left connected components of $S(X).$
\epro

\bpr left for the reader \epr


We also have (for algebraic homology):

\bthm\label{H0}
$H_0(X)=\Z^{r-1},$
where $r$ is the number of left orbits of $X$.
\ethm

\bpr
$H_0(X)$ is the quotient of the group
$$Z_0(X)=\left\{\sum_{x\in X} c_x (x): \sum_{x\in X} c_x=0\right\}$$
by its subgroup $B_0$ generated by $(x)-(yx).$
Let $x_0,...,x_{r-1}$ be representatives of all left orbits in $X$.
Let $\Z^{r-1}=\la e_1,...,e_{r-1}\ra$. We have $\phi: \Z^{r-1}\to Z_0(X)$ sending $e_i$ to $(x_i)-(x_0)$.
Since each generator of $B_0(X,\star)$ is a linear combination of elements of a single orbit, $\phi$ descends to an embedding $\psi:\Z^{r-1}\hookrightarrow H_0(X).$ It is easy to see that it is onto.
\epr

%
\subsection{Further properties of homology}
%

\blem\label{F_prop_lem} If $\star_1$ and $\star_2$ are mutually distributive then\\
(1) $(a_0\star_1 b)\star_2 (a_1\star_1 b)\star_2 ...\star_2 (a_n\star_1 b)=
(a_0\star_2 a_1\star_2 ...\star_2 a_n)\star_1 b.$\\
(2) $(a_0\star_1 ...\star_1 a_n)\star_2(a_k\star_1 ...\star_1 a_n)=
a_0\star_1 a_1\star_1 ...\star_1 a_{k-1}\star_{12} a_k\star_1 ... \star_1 a_n.$\\
(3) $(a_k\star_1 ... \star_1 a_n)\star_2 (a_{k+1}\star_1 ... \star_1 a_n)\star_2 ... (a_{n-1} \star_1 a_n) \star_2 a_n =a_k\star_{12} a_{k+1}\star_{12} ... \star_{12} a_n.$
\elem

\bpr (1) and (2) are by induction on $n.$\\
(3) By application of (1) to the left side of (3) and by induction on $n-k.$
\epr

For a binary product $\star$ on $X$ consider the function $F^\star: C_*(X)\to C_*(X)$
sending $(x_0,...,x_n)$ to
$$(x_0\star x_1\star ... \star x_n,...,x_{n-2}\star x_{n-1}\star x_n,x_{n-1}x_n,x_n).$$

Denote the composition $\star_1\star_2$ by $\star_{12},$ $\star_1\star_2\star_3$ by $\star_{123}$, and so on.

\bpro\label{F-prop}
(1) If $\star_1$ and $\star_2$ are mutually distributive and $\star_2$ is self-distributive then $\star_{12}$ is self-distributive and $F^{\star_1}$ is a chain homomorphism from $C_*(X,\star_{12})$ to
$C_*(X,\star_2)\\$
(2) If $\star_1,\star_2,\star_3$ are mutually distributive and $\star_3$ is self-distributive then
$F^{\star_{12}}= F^{\star_2}F^{\star_1}:C_*(X,\star_{123})\to C_*(X,\star_3)$
\epro

\bpr (1) Self-distributivity of $\star_{12}$ follows by direct computation.

Denote the $k$-th partial differentials on $C_n(X,\star_1)$ and $C_n(X,\star_{12})$
by $\p_{n,k}^1$ and $\p_{n,k}^{12}$, respectively. It is enough to prove that
\be{equation}\label{Fchain} F^{\star_1}\p_{n,k}^{12}(x_0,...,x_n)=\p_{n,k}^2F^{\star_1}(x_0,...,x_n),
\en{equation}
for all $k$ and $(x_0,...,x_n).$
The left hand side of (\ref{Fchain}) is
$$F^{\star_1}(x_0\star_{12}x_k,...,x_{k-1}\star_{12}x_k,x_{k+1},...,x_n)=
(w_0,...,w_{n-1}),$$
 where
$$w_i=(x_i\star_{12} x_k)\star_1 ... \star_1 (x_{k-1}\star_{12} x_k)\star_1 x_{k+1}\star_1... \star_1 x_n=$$
$$x_i\star_1 ... \star_1 x_{k-1}\star_{12} x_k\star_1 x_{k+1}\star_1... \star_1 x_n,$$
by Lemma \ref{F_prop_lem}(1) for $i<k$ and
$$w_i=x_{i+1}\star_1... \star_1 x_n,$$
for $k\leq i\leq n-1$.
The right side of (\ref{Fchain}) is
$$\p_{n,k}^2(x_0\star_1 ...\star_1 x_n,...,x_{n-1}\star_1 x_n,x_n)=(v_0,...,v_{n-1}),$$
 where
$$v_i=(x_i\star_1...\star_1 x_n)\star_2 (x_k\star_1...\star_1 x_n)=$$
$$x_i\star_1 ... \star_1 x_{k-1}\star_{12} x_k\star_1 x_{k+1}\star_1... \star_1 x_n=w_i,$$
by Lemma \ref{F_prop_lem}(3)
for $i<k$ and
$$v_i=x_{i+1}\star_1... \star_1 x_n=w_i,$$
for $k\leq i\leq n-1$.\\
(2) follows from  Lemma \ref{F_prop_lem}(3).
\epr

\brem
If $\star_1$ is invertible then, by Proposition \ref{F-prop}, $F^{\star_1}$ yields an isomorphism $F^{\star_1}_*:H_*(X,\star_{12})\to H_*(X,\star_2)$, despite the fact that the shelves $(X,\star_{12})$, $(X,\star_2)$ are, in general, not isomorphic.
For example, if $X=\{0,1\}$, $f(x)=1-x$ and $x\star_i y=f(x),$ for both $i=1,2,$ then $\star_{12}$ is the identity product which is not isomorphic to $\star_1.$
\erem


\bthm\label{vanishing}
If one of the conditions holds:\\
(1) $x\to x\star y$ is a bijection on $X$ for some $y$, or\\
(2) there is $y\in X$ such that $y\star x=y$ for all $x\in X$,\\
then $H_n(X)=0$ for all $n.$
\ethm

\bpr
(1) Since the map $f(x)=x\star y: X\to X$ is a shelf homomorphism, it induces
a chain map $C_*(X)\to C_*(X)$ which is chain homotopic to the zero map by the homotopy map $(x_0,...,x_n)\to (x_0,...,x_n,y)$. Therefore
$Id=f_*(f^{-1})_*=0$ as endomorphisms of $H_*(X)$ and, hence, $H_*(X)$ vanishes.\\
(2) If $y\star x=y$ for all $x\in X$ then the identity map on $C_*(X)$ is homotopic to the zero map via the chain homotopy
$C_n(X)\to C_{n+1}(X)$ sending $(x_0,...,x_n)$ to $(y,x_0,...,x_n).$
\epr

\bthm\label{xy=y}
If $X$ is a shelf with $x\star y=y$, c.f. Example \ref{shelf-ex}(2), then
$H_n(X)$ is a free abelian group with a basis given by $(x_0,x_1,...,x_n)-(x_1,x_1,x_2...,x_n)$
for all $x_0,..,x_n\in X,$ $x_0\ne x_1,$ $n>0.$ (In particular, $H_n(X)=\Z^{(|X|-1)|X|^n}$ for all $n\geq 0$ if $X$ is finite.)
\ethm

\bpr
Consider the map $s:C_n(X)\to C_{n+1}(X)$ sending $(x_0,...,x_n)$ to $(x_0,x_0,x_1,...,x_n).$ Note that $D_n=s(C_{n-1}(X))\subset C_n(X)$ defines a subcomplex $D_*$ of $C_*(X).$
Since $ds+sd$ is the identity on $D_n$, the identity map is homotopic to the zero map and, hence, this subcomplex is acyclic. Hence, the projection $\pi: C_*(X)\to C_*(X)/D_*$ induces an isomorphism $\pi_*: H_*(X)=H_*(C_*(X))\to H_*(C_*(X)/D_*).$
Since the differential on $C_*(X)/D_*$ vanishes, $H_n(X)=C_n(X)/D_n.$ Therefore, the elements $(x_0,x_1,...,x_n)-(x_1,x_1,x_2...,x_n)$ for $x_0\ne x_1$ are cycles in $C_n(X)$
whose images under $\pi_*$ form a basis of $H_*(C_*(X)/D_*).$ This implies the statement.
\epr


\subsection{Strong retracts}

\bthm
If $A$ is a strong retract of $X$ (c.f. Proposition \ref{s_retract}) then
$$H_n(X)=H_n(A)\oplus\Z(X\setminus A)\otimes \bigoplus_{k=0}^{n-1} H_k(A)\otimes \Z X^{n-k-1}.$$
\ethm

\bpr
Let $V$ be the subgroup of $\Z X$ generated by the elements $x-r(x)$ for $x\in X\setminus A.$ Let
$$C_{n,k}=\Z A^{n+1-k}\otimes V \otimes \Z X^{k-1}$$
be considered as a subgroup of $\Z X^{n+1-k}\otimes \Z X \otimes \Z X^{k-1}=\Z X^{n+1}=C_n(X),$ for $1\leq k\leq n+1,$ and let
$$C_{n,0}=\Z A^{n+1},\quad C_{n,k}=0\ \text{for $k>n+1$.}$$
If we assume additionally that $C_{-1,k}$ is $\Z$ for $k=0$ and it vanishes for $k>0$ then the chain complex $C_*(X)$ decomposes into a direct sum of subcomplexes
$$C_*(X)=\oplus_{k=0}^\infty C_{*,k}$$
and that the chain complex $C_{*,k}$ is isomorphic to $C_{*-k}(A)\otimes V\otimes \Z X^{k-1}$, for $k\geq 1,$ and to $C_*(A),$ for $k=0.$ (We assume here that $C_n(A)=0$ for $n<-1.$)
Replacing $k$ by $n-k$ implies the statement.
\epr

Denote the rank of an abelian group $A$ by $rk\, A.$

\bcor\label{cor_retract}
(1) If $A$ is a strong retract of a finite $X$, then
$$rk\, H_n(X)=rk\, H_n(A)+|X-A|\cdot \sum_{k=0}^{n-1} rk\, H_k(A)\cdot |X|^{n-k-1}.$$
(2) In particular, if $X$ is a shelf with $x\star y=g(y)$ for $g: X\to X$ such that $g^2=g$, c.f. Example \ref{shelf-ex}(2), then
$H_n(X)=\Z^{(r-1)|X|^{n}}$ for every $n$, where $r=|g(X)|.$\\
\ecor

More generally, Corollary \ref{cor_retract}(1) implies that if $rk\, H_k(A)=(c-1)|A|^k,$ where $c$ is the number of left orbits in $X$ for every $k$ then $rk\, H_n(X)=(c-1)|X|^n.$

\subsection{Computations and Conjectures}

Computer computations yield the following types of homology groups
for $3$-element shelves.

\begin{enumerate}
\item $H_n(X_r)=\Z^{(r-1)3^n},$ where $r=1,2,3$ is the number of left orbits, and\\
\item $H_n(X')=\begin{cases} \Z & \text{for $n=0$}\\ \Z^{2\cdot 3^{n-1}} & \text{for $n>0.$}\\ \end{cases}$
\end{enumerate}

The last type is realized for example by $X=\{1,2,3\}$ with multiplication
$$x\star y=\begin{cases} y & \text{if $y\geq x$}\\ 1 & \text{otherwise.}\end{cases}$$

Computer computations for $4$-element shelves yield the following types of the free parts of the homology groups:
\begin{enumerate}
\item $rk\, H_n(X_r)=(r-1)4^n,$ where $r=1,2,3,4$ is the number of left orbits, and, the ``exceptional" ones:
\item $rk\, H_0(X_1')=0$, $rk\, H_n(X_1')=4^{n-1}$ for $n>0.$
\item $rk\, H_0(X_2')=1$, $rk\, H_n(X_2')=2\cdot 4^{n-1}$ for $n>0,$
\item $rk\, H_0(X_3')=1$, $rk\, H_n(X_3')=3\cdot 4^{n-1}$ for $n>0,$
\item $rk\, H_0(X_4')=1$, $rk\, H_1(X_4')=3,$ $rk\, H_n(X_4')=13\cdot 4^{n-2}$ for $n\geq 2,$
\item $rk\, H_0(X_5')=2,$ $rk\, H_n(X_5')=7\cdot 4^{n-1}$ for $n>0.$
\end{enumerate}

On the basis of these computations we propose:

\bcon For every finite shelf $X$,\\
 $rk\, H_{n+1}(X)=|X|\cdot rk\, H_{n}(X)$ for $n\geq |X|-2.$
\econ

Furthermore our computations suggest:

\bcon For shelves in Example \ref{shelf-ex}(4), $H_n(X)$ is a group of rank
$$|X|^{n-1}(2+(|X|+1)(r-2))$$ for $n\geq 1,$ where $r$ is the number of left orbits in $X.$
\econ
This conjecture was proven for $n=1$ in \cite{PP}.

The above computations show that the map $\pi_*: H_*(X)\to H_*(O)$ induced by the homomorphism $\pi: X\to O$ of Remark \ref{r_orb} is (in general) neither $1$-$1$ nor onto.

Finally, it is worth noting that \cite{CPP,Pr} show two examples of $4$-element shelves which
have torsion in their first homology. Both of them are of the type described in Example \ref{shelf-ex}(4).

%
\section{Multi-shelf, Rack, and Quandle Homologies}
\label{s_multi}
%

Consider a multi-shelf $(X,\star_1,...,\star_n)$ (i.e. a set with mutually distributive products). Let $d_*^1,...,d_*^n$ be the corresponding induced differentials on $C_*(X).$

\blem
$d_{n-1}^k d_n^l=-d_{n-1}^l d_n^k$ for any $k,l\in \{1,...,n\}.$
\elem

The proof is a modification of that of Lemma \ref{dd=0}:
Let $$d_{n,i}^k(x_0,...,x_{n+1})=(x_0\star_k x_i,...,x_{i-1}\star_k x_i,x_{i+1},...,x_n).$$
Then $d_n^k=\sum_{i=0}^n (-1)^i d_{n,i}^k$ and it is easy to see that
$d_{n-1,i+1}^kd_{n,j}^l=d_{n-1,j}^ld_{n,i}^k,$ for $i\geq j$.
Hence, the $(n+1)n$ summands of $d_{n-1}^kd_n^l(x_0, ...,x_n)$ can be matched with coinciding terms of $d_{n-1}^ld_n^k(x_0, ...,x_n)$ of opposite signs.
\qed

Let $d_n=\sum_{k=1}^n c_k d_n^k :C_n(X)\to C_{n-1}(X)$ for some fixed $c_1,...,c_k\in \Z.$
By the above lemma, $d_{n-1}d_n=0$ for all $n$ and, hence, $(C_*(X),d_*)$ is a chain complex.
We denote its homology by $H_*(X, \sum_{k=1}^n c_k \star_k)$ and call it
the {\em multi-shelf homology}.

\bexa The identity product, $x\star_0 y=x$ is mutually distributive with every shelf product $\star$. The groups $H_*(X,\star - \star_0)$ are the rack homology groups introduced in \cite{FRS,CJKS2}, c.f. \cite{EG,LN}.
(This homology was originally defined for racks only, hence its name, even though the definition makes sense for all shelves.)
\eexa

\brem It is important to note however that by the standard definition of rack homology, $H_n(X, \star - \star_0)$ is the $n+1$-st rack homology of $X$. The connections between the shelf homology of $X$ and the simplicial complex of $X$ as well as Theorem \ref{H0} justify however, in our opinion, our choice of grading of the shelf homology.
\erem

By Proposition \ref{F-prop} we have

\bcor For every invertible $\star$ which is mutually distributive with multi-shelf
products $\star_1,...,\star_k$ on $X$, $F^{\star}$ defines an isomorphism between
$H_*(X, \sum_{k=1}^n c_k \star_k)$ and $H_*(X, \sum_{k=1}^n c_k \star\star_k).$
In particular, we have an isomorphism of rack homologies:
$H_*(X, \star - \star_0)\simeq H_*(X, \star^{-1} - \star_0),$
a property proved first by S. Kamada.
\ecor

Let $D_n(X)\subset C_n(X)$ be generated by chains $(x_0,...,x_n)$ such that $x_i=x_{i+1}$ for some $i.$ If a rack $X$ is a spindle (that is $x\star x=x$ for all $x\in X$) then $D_*(X)$ is a sub-chain complex of $C_*(X)$ of {\em degenerate chains}. Therefore one can consider the homology of $C_*(X)/D_*(X).$
The $n$-th homology of $C_*(X)/D_*(X)$ with the differential $\star - \star_0$ is called the {\em $n+1$-st quandle homology} of $(X,\star),$\cite{CKS,LN,Mo,NP1,NP2}.
(It is defined for all spindles despite the term ``quandle" in its name.) Applications of quandle homology to topology are discussed in \cite{Ca,CEGS, CJKLS,CJKS1,CJKS2,I,J,M,N,Za}.

Here is another interesting example of the multi-shelf homology:

\be{con} The Boolean multi-shelf, c.f. Example \ref{mshelf-ex}(2), that is $X=2^\Omega$ with
$$x_1*_0 x_2 =x_1,\quad x_1*_1 x_2=x_1\cap x_2,\quad x_1*_2 x_2=x_1\cup x_2$$
has homology groups of the following ranks:
$rk\, H_n(X,a_0\star_0 + a_1\star_1 +a_2\star_2)=$
$$\hspace*{.5in} \be{cases}
2^{|\Omega|\cdot (n+1)} & \text{for}\ (a_0,a_1,a_2)=(0,0,0)\\
|\Omega|\cdot 2^{n} & \text{for}\ (a_0,a_1,a_2)=c(1,-1,-1),\ c\in \Z\setminus \{0\}\\
1 & \text{for $a_0,a_1,a_2$ such that}\ a_0+a_1+a_2=0, (a_0,a_1,a_2)\ne(0,0,0) \\
0 & \text{otherwise,}\\
\en{cases}$$
for $n>0$ and
$$\hspace*{.5in} \be{cases}
2^{|\Omega|}-1 & \text{for}\ (a_0,a_1,a_2)=(0,0,0)\\
|\Omega| & \text{for}\ (a_0,a_1,a_2)=c(1,-1,-1),\ c\in \Z\setminus \{0\}\\
0 & \text{otherwise,}\\
\en{cases}$$
for $n=0.$
\en{con}

More generally we propose:

\be{con}
For every set $\star_1,...,\star_N$ of mutually distributive products,
there is a finite union of hyperplanes $H$ in $\Z^N$ and a sequence
$r_1,r_2,...\in \Z_{\geq 0}$
such that $rk\, H_n(X, a_1\star_1+ ... +a_N\star_N)=r_n$ for all $(a_1,...,a_N)\in \Z^N\setminus H$ and for all $n.$
\en{con}

In other words the above conjecture states that $rk\, H_n(X, a_1\star_1+ ... +a_N\star_N)$ is independent of $(a_1,...,a_N)$ with a codimension $\geq 1$ of possible exceptions.

\end{document}